\documentclass[a4paper]{jpconf}
\usepackage{amsmath,amssymb,graphicx}
\begin{document}
\title{MLD Relations of Pisot Substitution Tilings}

\author{Franz G\"ahler}

\address{Faculty of Mathematics, University of Bielefeld, D-33615 Bielefeld, 
Germany}

\ead{gaehler@math.uni-bielefeld.de}

\begin{abstract}
We consider 1-dimensional, unimodular Pisot substitution tilings with 
three intervals, and discuss conditions under which pairs of such tilings 
are locally isomorhphic (LI), or mutually locally derivable (MDL). For this
purpose, we regard the substitutions as homomorphisms of the underlying
free group with three generators. Then, if two substitutions are conjugated
by an inner automorphism of the free group, the two tilings are LI,
and a conjugating outer automorphism between two substitutions can often 
be used to prove that the two tilings are MLD. We present several examples
illustrating the different phenomena that can occur in this context. In 
particular, we show how two substitution tilings can be MLD even if their 
substitution matrices are not equal, but only conjugate in $GL(n,\mathbb{Z})$.
We also illustrate how the (in our case fractal) windows of MLD tilings 
can be reconstructed from each other, and discuss how the conjugating group
automorphism affects the substitution generating the window boundaries.
\end{abstract}

\section{Introduction}

In this article, we consider substitutions $\sigma$ on an alphabet of 
three letters, whose abelianisation matrix (substitution matrix) $M$ is 
primitive and unimodular, has irreducible characteristic polynomial, 
and a leading (Perron-Frobenius, PF) eigenvalue which is a Pisot number. 
The words generated by such a substitution can be regarded as elements 
of a free group with three generators, and automorphisms of the free 
group give rise to transformations of words. Alternatively, we can work 
with a geometric realisation of the substitution, by letting it act on 
three intervals, whose lengths are chosen proportional to the components 
of the left eigenvector associated with the leading PF-eigenvalue
$\lambda$ of $M$.  Each tile is then substituted with a sequence of
tiles, whose total length is equal to $\lambda$ times the original
length. Such a geometric realisation generates a tiling of the line,
instead of a sequence of symbols, or a word in a free group.

The matrix $M$ represents a linear mapping $A$ of $\mathbb{R}^3$,
expressed with respect to some basis $\{b_i\}$. As $M$ is unimodular,
this mapping is an automorphism of the lattice $L$ generated by
this basis.  We choose the geometry of $L$ such that the expanding and
contracting eigenspaces of $A$ are perpendicular to each other, so
that $A$ commutes with the orthogonal projections on these
eigenspaces. This can be realised as follows. After appropriate
rescaling, the tile lengths, being components of the PF-eigenvector,
are contained in the algebraic field $Q(\lambda)$, and so are all
coordinates of lattice points in the expanding eigenspace of $V$ of $A$. 
The corresponding coordinates in the contracting eigenspace $W$ can
be chosen as the $d-1$ Galois conjugates of the coordinate in $V$. 
We then have a cut-and-project scheme (CPS) defined by the lattice $L$, 
and the eigenspaces $V$ and $W$ of $A$:
\begin{equation}
\begin{matrix}
V\cong\mathbb{R}&\xleftarrow{\pi_1}&\mathbb{R}^3&\xrightarrow{\pi_2}&
  W\cong\mathbb{R}^2 \\
\cup      &                  &\cup        &                   &\cup    \\
\Lambda   &                  &L           &                   &\Omega
\end{matrix}
\label{gae_cps}
\end{equation}
One of the formulations of the Pisot conjecture states that the
vertex set $\Lambda$ of a Pisot substitution tiling always is a model
set, which means that there exists a window set $\Omega\subset W$
which is the closure of its interior, and which has boundary of
measure zero, such that $\Lambda = \{\pi_1(x)\,|\, x\in L,
\pi_2(x)\in\Omega\}$. Similarly, the subsets of the left end points of
all tiles of a given type in a Pisot substitution tiling are model
sets, too, with appropriate subwindows $\Omega_i$. For all examples
considered below, the Pisot conjecture can be shown to hold, even
though a proof for the general case is still missing. For a more
detailed description of Pisot substitution tilings and their associated
CPS, we refer to \cite{gae_frettl}. In particular, we note that one
can derive also a dual, contractive substitution acting on window sets 
in $W$, whose fixed points are the subwindows $\Omega_i$. We remark
that the CPS (1) also admits a canonical projection tiling, whose window
$\Omega$ is the image under the projection $\pi_2$ of the parallelepiped 
spanned by the basis $\{b_i\}$ of $L$. The subwindow for tile $i$, whose 
length is equal to the length of $\pi_1b_i$, is simply the parallelogram 
spanned by the vectors $\pi_2b_j$ and  $\pi_2b_k$, where $j$ and $k$ are 
the two indices different from $i$. This canonical projection tiling is 
not substitutional in general, but its subwindows can serve as convenient
starting points for the dual substitution determining the windows. Acting
on the canonical windows as seeds, the dual substitution overlap free.

\section{LI and MLD Relations}

The CPS (\ref{gae_cps}) does not specify the window $\Omega$ yet. 
Substitutions having the same abelianisation matrix $M$ (but differ 
in the order of the letters within a substituted word) give rise to 
the same CPS, but will have different windows in general. As we shall 
see below, even substitutions with different abelianisation matrices 
may belong to a common CPS.

In the following, we shall study relations between certain substitution 
tilings belonging to a common CPS. For this, besides the geometric 
realisation of a substitution tiling it is also useful to consider the 
substitution action on the underlying free group with three generators. 
In our examples, the substitution acts with a group automorphism.
If for two substitutions $\sigma_1$ and $\sigma_2$ there exists a fixed 
word $w$ in the group, such that $\sigma_1(g)=w^{-1}\sigma_2(g)w$ for 
every generator $g$ of the group, then the two substitutions produce 
tilings wich are locally isomorphic (LI), meaning that all their finite 
subpatterns are the same. This can be seen as follows. One first observes 
that there exists some power of $\sigma_1$, such that $\sigma_1^k$ has a 
bi-infinite fixed point, and that $\sigma_1^k$ and $\sigma_2^k$ are 
still conjugate in the same way, with a (longer) word $w'$. In a second 
step, one can then show that the fixed point of $\sigma_1^k$ is also a 
fixed point of $\sigma_2^k$, which implies that the two substitutions 
generate the same tilings.

A more delicate relation is mutual local derivability (MLD) \cite{gae_mld}. 
Two tilings are MLD, if one can be reconstructed from the other in a 
{\it local} way, and vice versa. For this to work, the two tilings must
first be brought to the appropriate relative scale and position. A good
starting point is to consider two tilings belonging to a common CPS. 
In fact, two (model set) tilings are MLD if and only if the window
of one can be constructed by finite unions and intersections of lattice
translates of the window of the other, and vice versa. Looking at the
windows can suggest an MLD relation, but for proving such a relation
it is very helpful if one substitution can be written as a conjugate
of the other, $\sigma_1=\rho^{-1}\circ\sigma_2\circ\rho$, where
$\rho$ is an outer automorphism of the free group (an inner automorphism
would lead to an LI relation). Such an automorphism will make the
transformation of one tiling into the other explicit. 

If a substitution $\sigma$ acts invertibly on the underlying free group, 
the boundaries of its windows are generated by a substitution, too.
This boundary substitution is given by $\sigma_b=\tilde{\sigma}^{-1}$, 
which is the inverse of $\sigma$, read backwards \cite{gae_ei,gae_arnoux2}.
As a seed for the iteration, it is again convenient to take the windows of 
the canonical projection tiling belonging to the same CPS. For instance, 
the canonical window for tile $a$ is bounded by the closed path consisting 
of the four consecutive segments $\pi_2b_2$, $\pi_2b_3$,$-\pi_2b_2$, and 
$-\pi_2b_3$. This path is represented by the word $bcb^{-1}c^{-1}$ in the 
free group, on which the boundary substitution $\sigma_b$ acts. In each 
step, the boundary path is transformed into one with more, but shorter 
segments, eventually converging to the fractal boundary of the final window.

If we now have two substitutions $\sigma$ and $\sigma'$ which are conjugated,
$\sigma'=\rho^{-1}\circ\sigma\circ\rho$, their boundary substitutions 
satisfy $\sigma_b'=\tilde\rho^{-1}\circ\sigma_b\circ\tilde\rho$. Iterating 
this, we find ${\sigma'}_b^n=\tilde\rho^{-1}\circ\sigma_b^n\circ\tilde\rho$.
In the limit of a fully fractalized window, the action of $\tilde\rho^{-1}$ 
can be neglected (it is local at the scale of the then infinitesimally 
small segments), and $\sigma'_b$ can be understood as $\sigma_b$ acting 
on the seed of $\sigma'_b$, transformed by $\tilde\rho$. As the dual
substitution acting on the canonical windows and their iterates is 
overlap free, the transformation $\tilde\rho$ commutes with the 
fractalization induced by $\sigma_b$, and it becomes manifest that the
windows of $\sigma$ are obtained from those of $\sigma'$ via the
transformation induced by $\tilde\rho$. In particular, the windows of the 
two substitutions have the same fractal structure.

In the following, different phenomena arising in this context are 
illustrated with a number of examples. As a short-hand notation, 
we write the action of a substitution $\sigma$ on a free group as the 
list of images of the generators, in our case a triple 
$[\sigma(a),\sigma(b),\sigma(c)]$.

\begin{figure}
\centerline{
\begin{minipage}{7cm}
\centerline{\includegraphics[height=6cm]{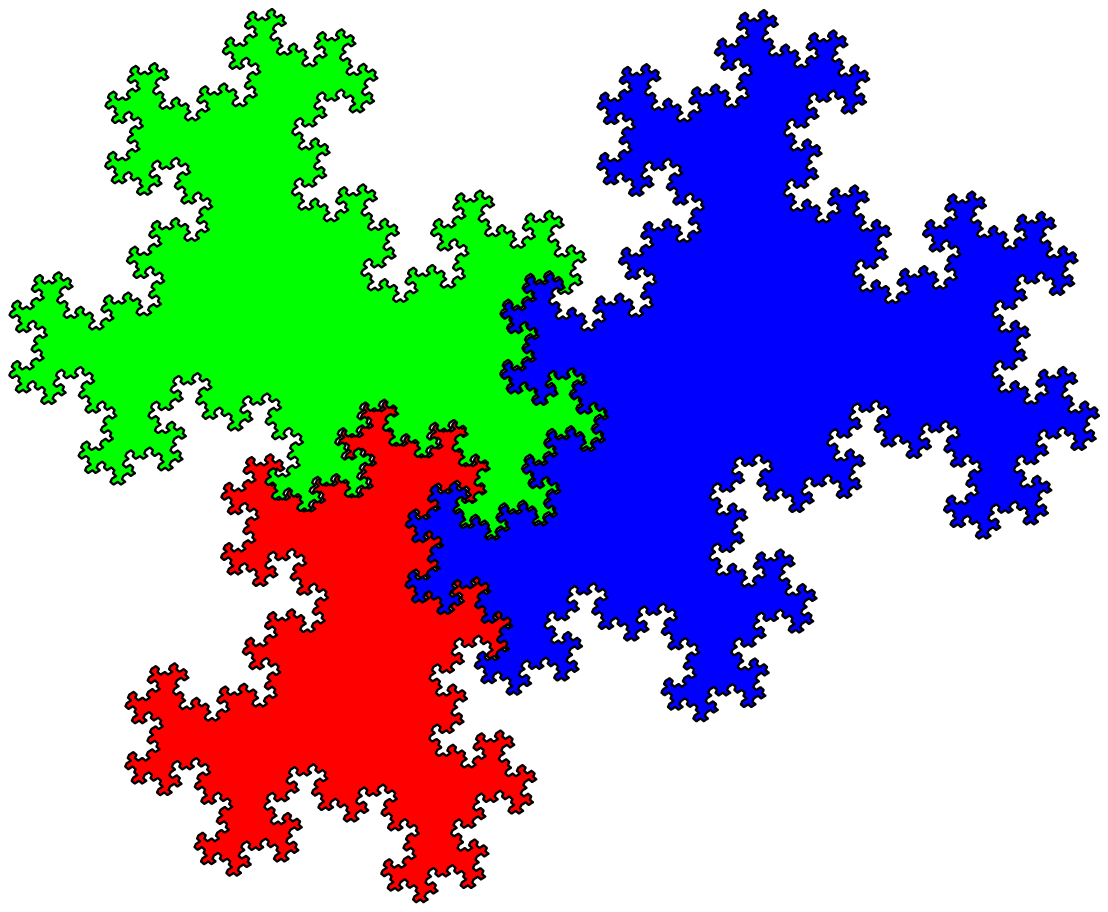}}
\caption{Windows of the substitution 
         $a\rightarrow cb$, $b\rightarrow c$, $c \rightarrow cab$.}
\end{minipage}\hspace{1cm}%
\begin{minipage}{7cm}
\centerline{\includegraphics[height=6cm]{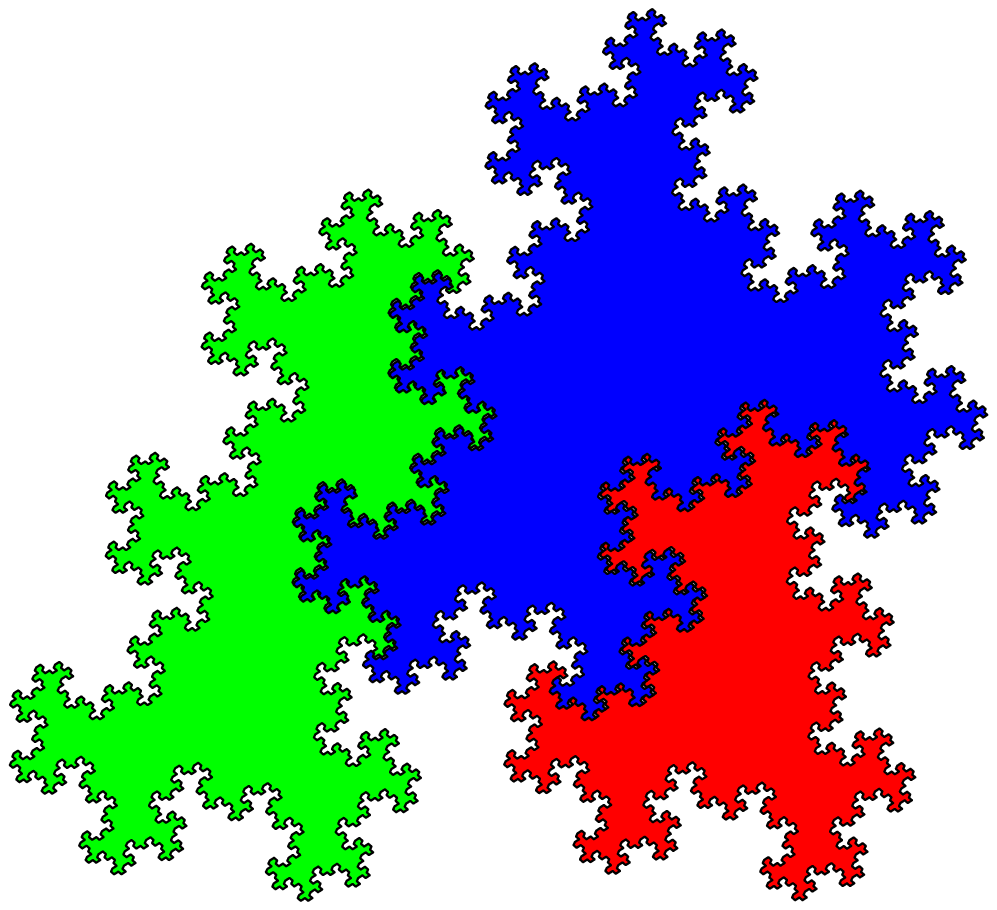}}
\caption{Windows of the substitution 
         $a\rightarrow bc$, $b\rightarrow c$, $c \rightarrow cba$.}
\end{minipage}
}
\end{figure}

\section{Examples}

As a first example, we consider the substitutions $\sigma_1=[cb,c,cab]$
and $\sigma'_1=[bc,c,cba]$, which have the same abelianisation matrix. 
These two substitutions are conjugate by the free group automorphism 
$\rho_1 = [bab^{-1},b,c]$, with inverse $\rho_1^{-1}=[b^{-1}ab,b,c]$. 
It is easily checked that indeed we have $\sigma_1=\rho_1^{-1}\circ
\sigma'_1\circ\rho_1$. In a word generated by $\sigma_1$, there is always 
a $b$ to the right of an $a$. $\rho_1$ eats up that $b$, and adds a $b$ 
to the left of the $a$ instead, effectively replacing $ab$ pairs by $ba$ 
pairs. $\rho_1^{-1}$ performs the opposite operation. This is obviously a 
local operation, no matter whether one works with words in a free group, 
with symbolic sequences, or with tilings. The LI classes of tilings 
generated by the two substitutions are MLD. The windows of the two
substitutions $\sigma_1$ and $\sigma'_1$ are shown in Figures 1 and 2,
respectively. The windows for the $a$, $b$, and $c$ tiles are in red 
(medium gray), green (light gray), and blue (dark gray). 
When transforming from Figure~1 to Figure~2, part of the $b$ tiles 
(green), namely those to the right of an $a$ tile, move to where the 
$a$ tiles were before. The subwindow of the $a$ tiles (red) thus becomes 
green, and a congruent copy of it is cut away from the original subwindow
of the $b$ tiles in green. The red subwindow of the $a$ tiles instead
moves to a different place, because the $a$ tiles are now to the right
of a $b$ tile. As $a$ and $b$ tiles have different lengths, the left
endpoint of the second tile of $ab$ and $ba$ pairs differs, and so 
the corresponding subwindows are at different places.

In the second example, we consider two substitutions with different
abelianisation matrices, $\sigma_2=[c,a,cab]$ and $\sigma'_2=[c,ca,cb]$.
Again, there is a conjugating automorphism $\rho_2=[a,a^{-1}b,c]$, with
inverse $\rho_2^{-1}=[a,ab,c]$, so that $\sigma_2=\rho_2^{-1}\circ
\sigma'_2\circ\rho_2$. Here, in words produced by $\sigma_2$, all
$b$ tiles are to the right of an $a$ tile. $\rho_2$ eats up the $a$ tile
to the left of a $b$ tile, effectively replacing all $ab$ pairs by just 
one $b$. Other $a$ tiles (not to the left of a $b$) and all $c$ tiles 
are left as they are. Conversely, $\rho_2^{-1}$ splits all $b$ in a 
$\sigma'_2$-word into $ab$ pairs. On the tiling level, this operation 
is local if and only if the length of an $ab$ pair of tiles in the 
$\sigma_2$-tiling is the same as the length of a $b$ tile in the 
$\sigma'_2$-tiling, whereas $a$ and $c$ tiles have the same length for 
both tilings. With appropriate global scalings, this is indeed the case. 
$\sigma_2=\rho_2^{-1}\circ\sigma'_2\circ\rho_2$ implies that the two 
abelianisation matrices are conjugate in $GL_3(\mathbb{Z})$. In fact, 
the two substitutions have the same CPS, with the same lattice $L$. 
The only difference is, that the linear mapping $A$ is expressed with 
respect to two different lattice bases, yielding two different matrix 
representations of $A$, and different tile lengths (which are the lengths 
of the projected basis vectors). It is therefore not surprising, that 
the length of tile $b$ in the $\sigma'_2$-tiling is the sum of the 
lengths of the two tiles $a$ and $b$ of the $\sigma_2$-tiling. The
windows of the substitutions $\sigma_2$ and $\sigma'_2$ are shown
in Figures~3 and 4, respectively, using the same coloring as for the
previous example. In Figure~3, part of the $a$ tiles (in red), namely
those to the left of a $b$ tile, become the new $b$ tiles in Figure~4 
(green), whereas the old $b$ tiles in Figure~3 (green) are discarded. 
MLD relations can therefore arise also if the two abelianisation matrices 
are not equal, but conjugate in $GL_3(\mathbb{Z})$, because the two 
substitutions then share a common CPS. We emphasise, however, that this 
relation is local only for the tilings with properly sized tiles. 
This pair of examples had been discussed in detail already in 
\cite{gae_arnoux}. $\sigma'_2$ is LI to the Rauzy or Tribonacci 
substitution.

\begin{figure}
\centerline{
\begin{minipage}{7cm}
\centerline{\includegraphics[height=6cm]{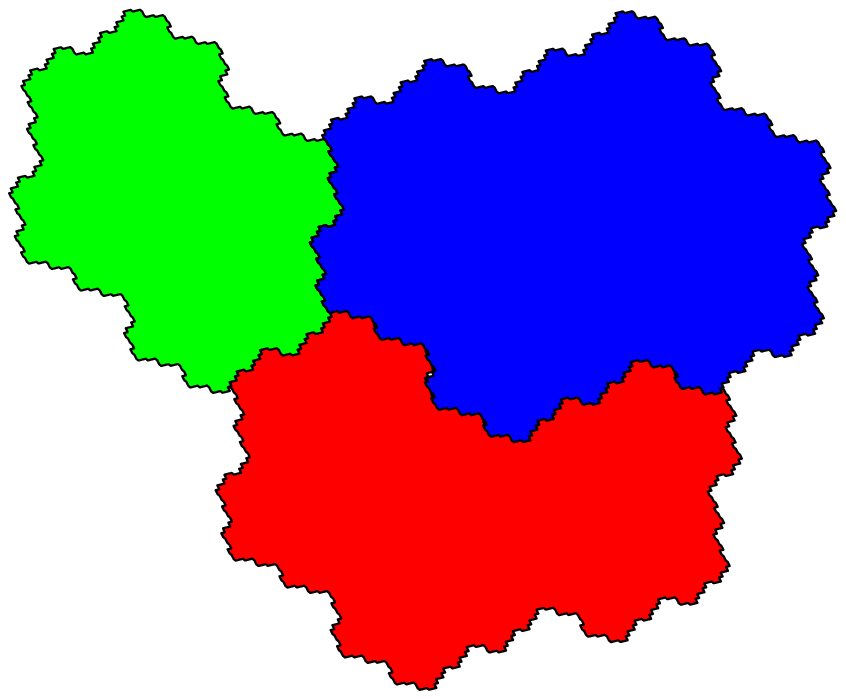}}
\caption{Windows of the substitution 
         $a\rightarrow c$, $b\rightarrow a$, $c \rightarrow cab$.}
\end{minipage}\hspace{1cm}%
\begin{minipage}{7cm}
\centerline{\includegraphics[height=6cm]{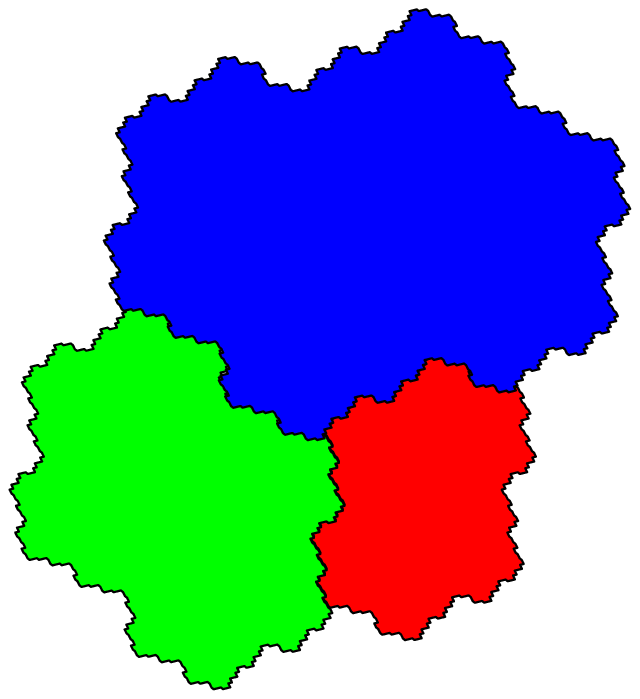}}
\caption{Windows of the substitution 
         $a\rightarrow c$, $b\rightarrow ca$, $c \rightarrow cb$.}
\end{minipage}
}
\end{figure}

Finally, as a third example, we consider a quartet of substitutions,
all with the same abelianisation matrix $M$. These substitutions are 
$\sigma_A=[ca,ab,cab]$, $\sigma_B=[ac,ab,abc]$, $\sigma_C=[ca,ba,bac]$, 
and $\sigma_D=[ac,ab,bac]$. $\sigma_A$ and $\sigma_D$ are conjugate in 
a way already seen in the first example: 
$\sigma_D=\rho_3^{-1}\circ\sigma_A\circ\rho_3$, where
$\rho_3=[a,b,a^{-1}ca]$ simply replaces $ac$ pairs by $ca$ pairs. 
In order to discuss the relations to the other substitutions,
we introduce the automorphisms $u_1=[c,a,ab]$ and $u_2=[c,a,ba]$. We 
then have $\sigma_A=u_1\circ u_1\circ u_2$, $\sigma_B=u_1\circ u_2\circ u_1$, 
and $\sigma_C=u_2\circ u_1\circ u_1$, so that $\sigma_A=u_1\circ\sigma_B
\circ u_1^{-1}$ and $\sigma_C=u_1^{-1}\circ\sigma_B\circ u_1$. The 
conjugating automorphism $u_1$ has an abelianisation matrix $U$ which 
commutes with the common abelianisation matrix $M$ of the substitutions, 
even though $U$ is non-trivial. This is possible because $M$ is equal to
the third power of $U$. Therefore, $M$ is conjugate to itself by some 
non-trivial mapping, which acts non-trivially on the lattice $L$, changing 
the scale of the tiling by the cubic root of the inflation factor $\lambda$ 
of the substitution (or its inverse). Consequently, in order to be MLD, the 
tilings produced by $\sigma_A$, $\sigma_B$ and $\sigma_C$ must be at relative
scales $\lambda^{-\frac13}$, $1$, and $\lambda^{\frac13}$, respectively. 
The situtation is in fact similar to the second example, where the 
substitutions share a common CPS, but different lattice bases of $L$ 
are used. Here, these different lattice bases still lead to the 
same abelianisation matrix $M$, but produce tiles of different sizes.
The windows of the substitutions $\sigma_A$, $\sigma_B$, $\sigma_C$,
and $\sigma_D$ are shown if Figures~5 to 8, respectively. We don't
discuss the transformations between them in detail here, but it is
quite obvious that each subwindow can be obtained as translate or
union of translates of subwindows of the other substitutions. 

\begin{figure}
\centerline{
\begin{minipage}{7cm}
\includegraphics[height=6cm]{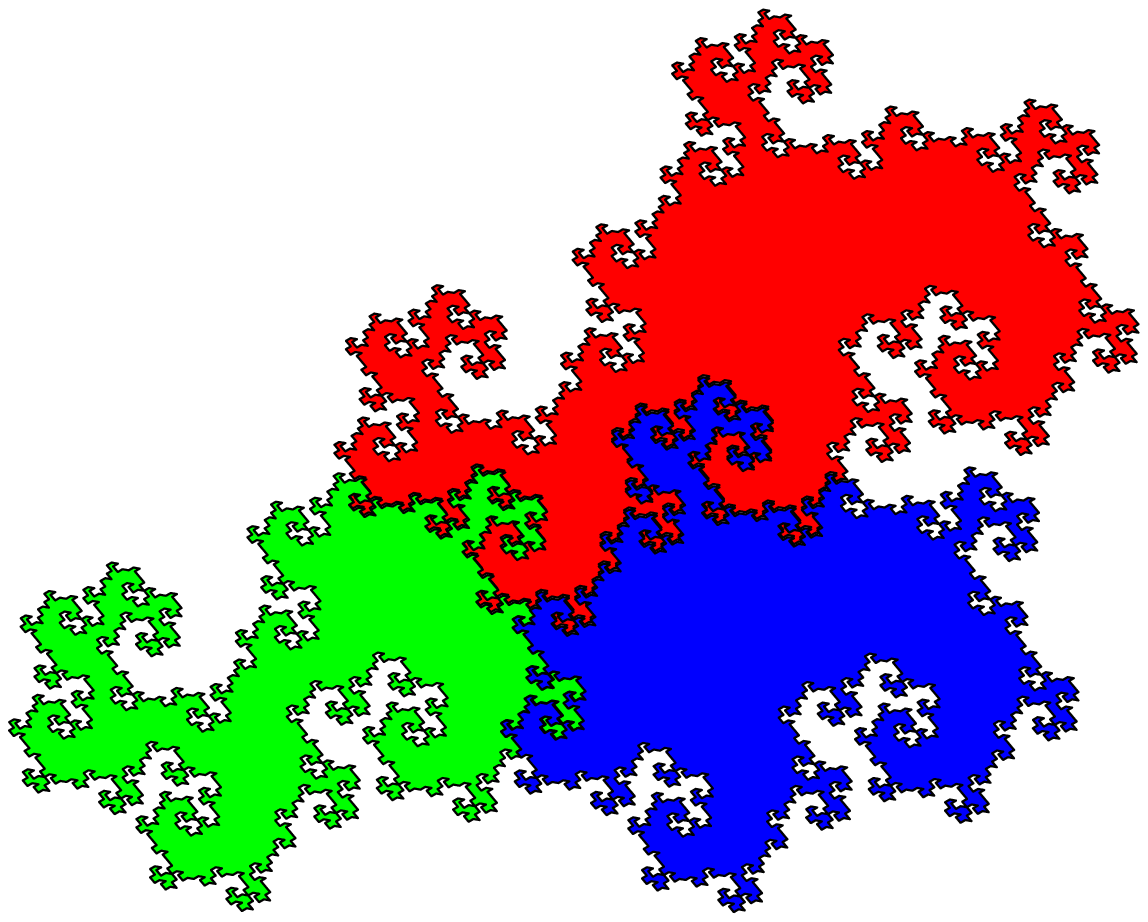}
\caption{Windows of the substitution 
         $a\rightarrow ca$, $b\rightarrow ab$, $c \rightarrow cab$.}
\end{minipage}\hspace{1cm}
\begin{minipage}{7cm}
\centerline{\includegraphics[height=6cm]{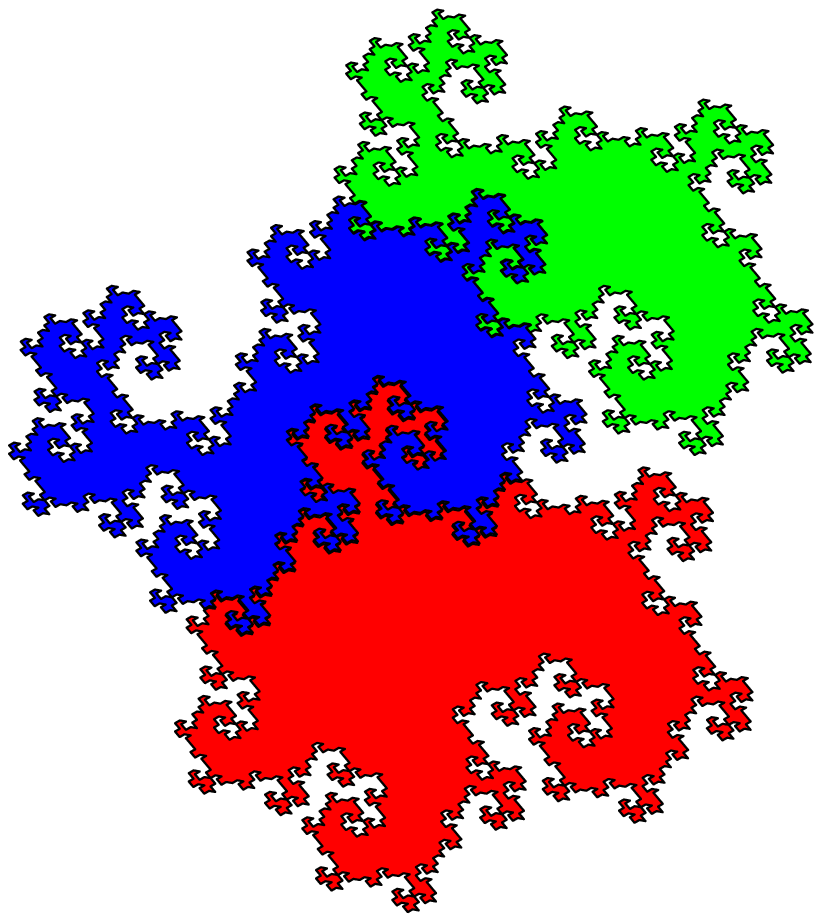}}
\caption{Windows of the substitution 
         $a\rightarrow ac$, $b\rightarrow ab$, $c \rightarrow abc$.}
\end{minipage}
}
\null\vspace{1cm}
\centerline{
\begin{minipage}{7cm}
\includegraphics[height=6cm]{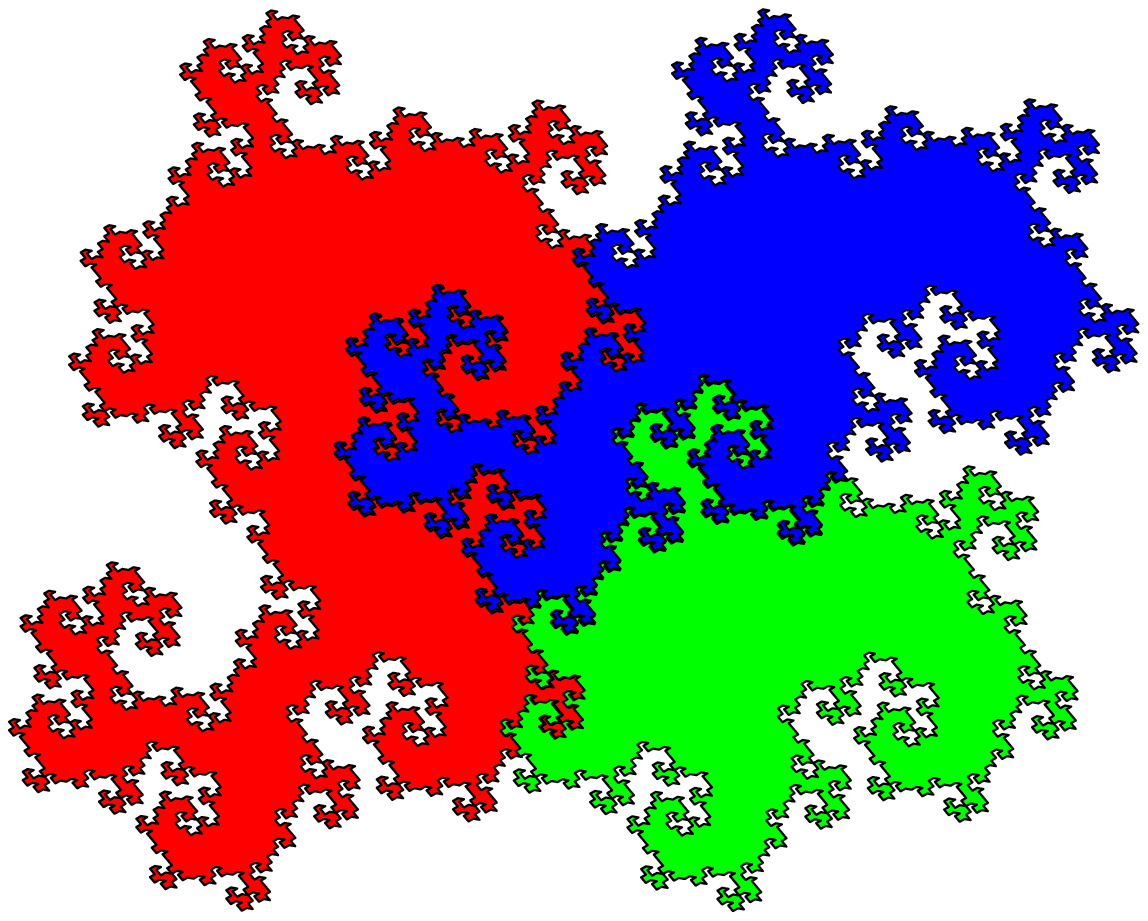}
\caption{Windows of the substitution 
         $a\rightarrow ca$, $b\rightarrow ba$, $c \rightarrow bac$.}
\end{minipage}\hspace{1cm}
\begin{minipage}{7cm}
\centerline{\includegraphics[height=6cm]{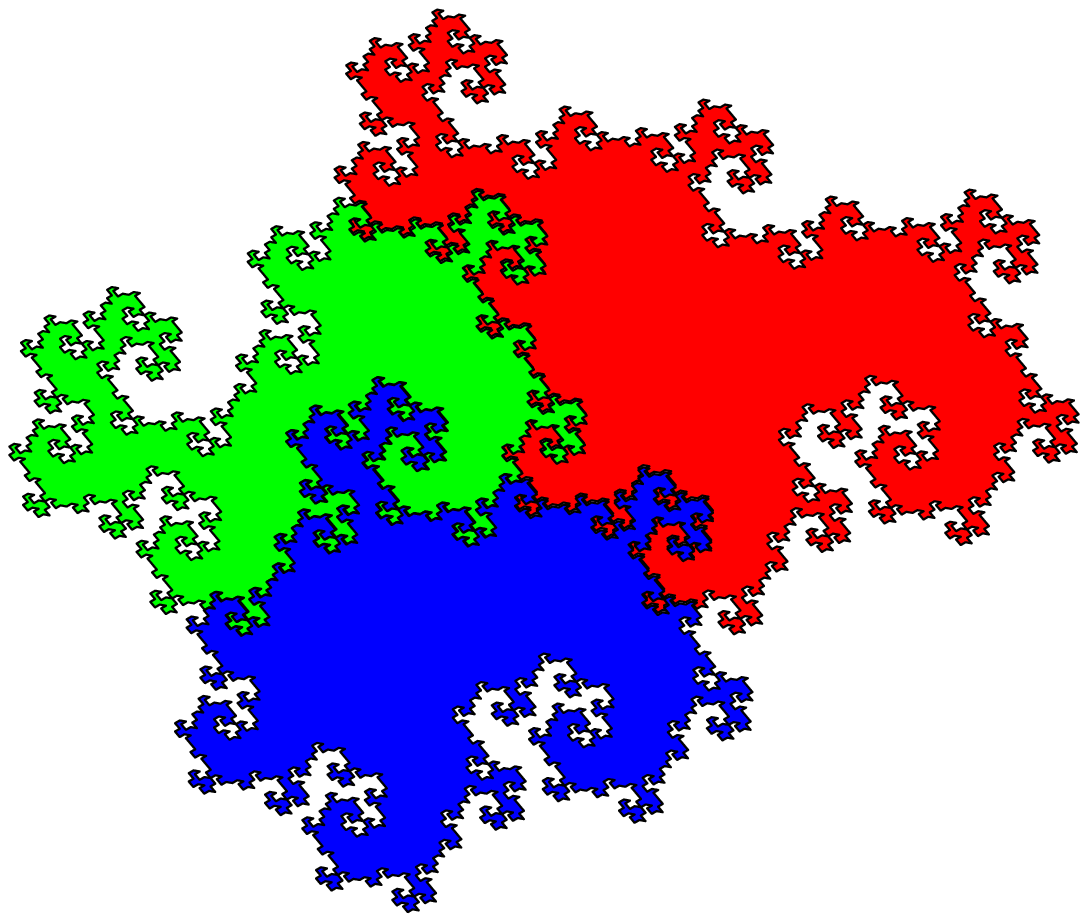}}
\caption{Windows of the substitution 
         $a\rightarrow ac$, $b\rightarrow ab$, $c \rightarrow bac$.}
\end{minipage}
}
\end{figure}

\section{Conclusions}

The CPS of a Pisot substitution tiling can accommodate many other
substitution tilings as well. Some of these are obtained by
permuting the letters in the substituted words, but there may be
others arising from the choice of a different basis for the lattice 
of the CPS, as we have seen in the second example. The more
complicated a substitution is, the richer is the set of substitution
tilings supported by its CPS. Some of the tilings sharing a common 
CPS can be related in interesting ways, however. In particular, there 
may local isomorphism or mutual local derivability relations, sometimes
in surprising ways. We have discussed some of the phenomena that may 
arise in this context, and have illustrated them with a number of 
examples. The key tool was to formulate the substitution as an
automorphism of an underlying free group, which allowed to find
the relations with algebraic methods, and to make the transformations
between the related tilings explicit.

\ack 
The author would like to thank Pierre Arnoux, Dirk Frett\-l\"oh, and 
Edmund Harriss for fruitful discussions.

\end{document}